\documentclass[12pt]{amsart}
\usepackage{amsmath,amsthm}

\def\R{{\mathbb R}}

\def\d{{\rm d}}
\def\p{\partial}

\def\l{\lambda}
\def\r{\rho}
\def\w{\omega}
\def\a{\alpha}

\def\T{{\mathcal T}}
\def\M{{\mathcal M}}

\def\o{\overline}

\def\WP{Weil-Petersson\ }

\theoremstyle{plain}
\newtheorem{thm}{Theorem}[section]
\newtheorem{lm}[thm]{Lemma}
\newtheorem{prop}[thm]{Proposition}
\newtheorem{cor}[thm]{Corollary}
\newtheorem{fact}[thm]{Fact}

\theoremstyle{definition}
\newtheorem{df}[thm]{Definition}

\numberwithin{equation}{section}

\title{Upper bound for the Weil-Petersson volumes}
\author{Draft by: Samuel Grushevsky}
\email{grushevs@math.harvard.edu}
\address{Mathematics Department, Harvard University
One Oxford Street, Cambridge, MA 02138, USA}
\thanks{Partially supported by NSF Graduate Research Fellowship}
\date{December 9, 1999}
\begin{document}
\setcounter{section}{-1}

\begin{abstract}
An explicit upper bound for the \WP volumes of punctured Riemann
surfaces is obtained  using the combinatorial integration scheme
from \cite{wpvol}. It is shown that for a fixed number of
punctures $n$ and for genus $g$ increasing,
$$\lim\limits_{g\to\infty, \ n\ {\rm fixed}}\frac{\ln {\rm vol}_{WP}
(\M_{g,n})}{g\ln g}\le 2,$$ while this limit is exactly equal to
two for $n=1$.
\end{abstract}
\maketitle

\section{Introduction}
After Wolpert in \cite{wolpert} computed the cohomology of the
moduli space of Riemann surfaces as a graded vector space, the
question of computing the cohomology ring structure (aka the
intersection theory) on the moduli arose. The problem has been
intensively studied since then. Witten's paper \cite{witten} is a
good source for available techniques and ideas. Witten's
conjecture, which later became Kontsevich's theorem
\cite{kontsevich}, shows that the intersection numbers satisfy a
certain KdV equation.

However, the problem of getting explicit numerical results still
remains, as the recursive computations become exceedingly complicated as
the genus and number of punctures grow. Carel Faber computed some 
low-genus intersection numbers in \cite{fabnum} and has obtained numerous
results in other papers.

Recently Zograf \cite{zograf} and Manin and Zograf \cite{mz} have
obtained quite explicit generating functions for the \WP volumes, and
computed the asymptotics of the volume growth for genus being fixed and
the number of punctures growing to infinity.

In this paper we use a completely different set of tools, namely the
decorated Teichm\"uller theory, to obtain an explicit asymptotic upper bound for
the \WP volumes for a fixed number of punctures and the genus growing to
infinity.

\section{Decorated Teichm\"uller Theory}
Let us recall the notations and relevant constructions. All of
these come from Penner's work \cite{wpvol}: this is just a brief
summary.

Let $\M_{g,n}$ denote the moduli space of Riemann surfaces of genus $g$ with
$n$ punctures --- it has complex dimension $3g-3+n$. Let $\w_{WP}$ denote
the two-form of the \WP scalar product on $\M_{g,n}$. It can be extended
to a closed current on the Deligne-Mumford compactification  $\o{
\M_{g,n}}$ of the moduli space. Taking its highest power produces a
volume form on $\M_{g,n}$, integrating which over $\o{\M_{g,n}}$ gives
the {\it \WP volume} ${\rm vol}_{WP}(\M_{g,n})$, which is the principal object
of our study. We will only be concerned with hyperbolic punctured
surfaces, i.e.~the case when $2g+n\ge 3$ and $n>0$.

An {\it ideal triangulation} of a punctured surface is a triangulation
of the surface with vertices only at punctures. We can straighten an
ideal triangulation so that every edge of it is a geodesic arc for the
hyperbolic metric on the surface. From Euler characteristics
considerations it follows that such a triangulation consists of 
$V:=4g-4+2n$ triangles, and has $N:=6g-6+3n$ edges.

The {\it decoration} of a punctured Riemann surface is an addition of a
horocycle around each puncture. More rigorously, on the uniformizing
hyperbolic plane we take a horocycle around a preimage of a puncture,
and consider its projection to the Riemann surface. For a decorated
Riemann surface with an ideal geodesic triangulation define the 
$\l$-length of arc $e$ of the triangulation to be $\l(e)=\sqrt{2e^\delta}$,
where $\delta$ is the (signed) hyperbolic distance from the point where
$e$ intersects the horocycle around one its end to the point of
intersection with the horocycle at the other end. It turns out 
(\cite{intscheme}, theorem 3.1) that for
a fixed triangulation the $\l$-lengths establish a homeomorphism of the
decorated Teichm\"uller space $\widetilde{\T_{g,n}}$ and $\R_+^{6g-
6+3n}$.

An {\it embedded graph} is a graph embedded in a Riemann surface.
Combinatorially it can be represented as a usual graph endowed with a
cyclic order of the edges around each vertex. Given an ideal
triangulation of a surface, taking its Poincar\'e dual produces a
trivalent embedded graph on the surface, which we denote $\Gamma$. The graph
$\Gamma$ has $N$ edges and $V$ vertices. We
define $\l(e)$, the {\it $\l$-length of an edge $e$ of the graph}, to be the 
$\l$-length of its Poincar\'e dual geodesic arc of the ideal triangulation.

\section{Moduli space description}
Our goal is to use Penner's description of the moduli space in the 
$\l$-length coordinates to obtain an explicit upper bound of the \WP volumes.

For any edge $e$ of the trivalent embedded graph $\Gamma$ let $f_i$ and
$g_i$ be the adjacent edges of the graph at $i$'s end of $e$. Then
define the associated simplicial coordinate to be
$$X_e:=\sum\limits_{i=1}^2\frac{\l(f_i)}{\l(e)\l(g_i)}+\frac{\l(g_i)}
{\l(e)\l(f_i)}-\frac{\l(e)}{\l(f_i)\l(g_i)}.$$
Further, let $\r_i$ be the sum of simplicial coordinates of the edges in
a path around puncture number $i$, where the path has to go to ``the
next edge to the left'' at each vertex in the cyclic order corresponding
to the embedding of the graph. Notice that if we take any edge and start
going to ``the next edge to the left'' from it, we will end up with a loop
around some puncture. Since you can ``go to the left'' in two directions,
i.e.~since it matters in which way you start going, we have
$$\r:=\sum\limits_{i=1}^n\r_i=2\sum_{i=1}^N X_{e_i}=2\sum
\frac{\l(e)}{\l(f)\l(g)},$$
where the sum is over all triples of edges having a common vertex,
including the possible renamings of $e$, $f$ and $g$.

In these notations Penner proves the following result:

\begin{fact}[\cite{wpvol}, 3.2.1 and 3.4.3]
\label{moduli}
Let $\w$ be a top-dimension differential form on $\M_{g,n}$. Then
\begin{equation}
\int_{\o\M_{g,n}}\w=\sum_{[\Gamma]}\frac{1}{{\rm Aut}\,
\Gamma}\int_{D(\Gamma)}\pi^*(\w),
\label{domain}
\end{equation}
where $\pi:\widetilde\T_{g,n}\to\M_{g,n}$ is the forgetful projection,
$[\Gamma]$ denotes the isomorphism class of an embedded trivalent graph
$\Gamma$, and $D(\Gamma)$ is the domain in the decorated Teichm\"uller
space given by
\begin{equation}
\label{domainD}
D(\Gamma)=\left\lbrace \r_i=1;\ X_e>0|\forall e\in\Gamma,\ \forall
i=1\ldots n\right\rbrace
\end{equation}
in the $\l$-coordinates corresponding to the embedded graph $\Gamma$.
\end{fact}

\section{\WP Volume form}
In further computations for simplicity we drop the $\l$'s
and simply write $e$ for $\l(e)$, if no confusion is possible. 

In $\l$-length coordinates on the moduli the \WP two-form is given by
(\cite{wpvol}, theorem A.2)
\begin{equation}
\label{twoform}
\w_{WP}=-2\sum\limits_{v\in\Gamma;\ e,f,g\ni v} \frac{\d e}{e}\wedge
\frac{\d f}{f}+\frac{\d f}{f}\wedge\frac{\d g}{g}+\frac{\d g}{g}\wedge
\frac{\d e}{e}.
\end{equation}

The \WP volume form is the $3g-3+n$'s external power of $w_{WP}$. In general,
letting $I$ be a multi-index, and denoting the exclusion of factors by a
hat, it would be the sum
\begin{equation}
\label{volform}
\w_{WP}^{\wedge(3g-3+n)}=\sum\limits_{|I|=n} a_I \prod
d\ln\l_1\wedge\ldots\wedge\widehat{d\ln\l_I}\wedge\ldots\wedge d\ln\l_N.
\end{equation}
\begin{prop}
\label{aI}
In the above notations, $|a_I|\le 2^N$.
\end{prop}
\begin{proof}
Use the expression (\ref{twoform}) for the two-form
to straightforwardly take an exterior power and compute the total number
of summands of one kind (with fixed $I$). Suppose the product contains
some $d\ln e$; then it must come in pair with one of the adjoining
edges --- let it be $f_i$ in the above notations, so there are four
choices. Then $g_i$ must come (if it is in the sum) with one of the two
edges at its other end, and then the other edge at its other end must
come with still another, and so on. Thus we have one factor of four, and
many factors of two. Each new factor of four appears if we encounter an
edge in $I$, so the total number of summands of one kind is at most $2^|I|
2^{3g-3+n}=2^{3g-3+2n}$. Recalling the $-2$
in the two-form, $|a_I|\le 2^{3g-3+2n}\cdot 2^{3g-3+n}=2^N$.
\end{proof}
For the case of one puncture Penner explicitly computes the \WP volume form to
be (\cite{wpvol}, theorem 6.1.2)
\begin{equation}
\label{allvarform}
\w_{WP}^{\wedge(3g-2)}=\pm 2^{4g-2}\sum\limits_{i=1}^N (-1)^i d\ln
\l_1\wedge\ldots\wedge\widehat{d\ln\l_i}\wedge\ldots\wedge d\ln\l_N.
\end{equation}
Thus we are dealing with a form which is singular when the $\l$
coordinates approach zero. However, this is not a problem:
\begin{prop}
\label{ge4}
In the domain of integration $D(\Gamma)$ (formula \ref{domain})
 for any edge $e\in\Gamma$ we have $\l(e)>4$ in $D$.
\end{prop}
\begin{proof}
We can construct at most two paths
around punctures ``going to the left'' including the edge $e$ --- a
path can be constructed by deciding at which end of $e$ we start
building it. Let one of these paths go through edges $f_1, e, f_2$, and
the other --- through $g_1,e,g_2$; let $\rho_i$ and $\rho_j$
be the sums of the simplicial coordinates of the edges in these paths.
Then using the formula for $\rho$ in terms of $\alpha$-lengths of
sectors (\cite{wpvol}, lemma 3.4.2) we see that
$$\rho_i>\frac{2}{e}\left(\frac{g_1}{f_1}+\frac{g_2}{f_2}\right)\quad{\rm
and}\quad\rho_j>\frac{2}{e}\left(\frac{f_1}{g_1}+\frac{f_2}{g_2}\right).$$
Since $\rho_i=\rho_j=1$ in the domain $D$, we get
$$2>\frac{2}{e}\left(\frac{f_1}{g_1}+\frac{g_1}{f_1}+\frac{f_2}{g_2}+
\frac{g_2}{f_2}\right)>\frac{8}{e},$$
so that $e>4$. If there were only one path going through $e$, i.e.~if
$i=j$, we would have $1=\rho_i>\frac{8}{e}$, and thus $e>8$, which is
even better.
\end{proof}

Thus in the domain of integration $D(\Gamma)$ the $\l$-lengths are
bounded below. However, $D$ can have limit points at infinities of $\l$-
lengths, and since the integral of $d\ln x$ does not converge at
infinity, we have to deal with this problem in great detail, using
the triangle inequality.

\section{Triangle Inequality}
The problem with the decorated Teichm\"uller theory is that the domain
of intgeration $D(\Gamma)$ cannot be simply described in $\l$-lengths.
Our success will come from the following observation:

\begin{thm}[Triangle inequality: \cite{intscheme}, lemma 5.2]
\label{le2thm}
Let $e$, $f$, and $g$ be three  edges of the graph $\Gamma$ having a common
vertex. Then in the domain of integration $D(\Gamma)$ the triangle inequality
between them holds:
\begin{equation}
\label{triangneq}
e\le f+g.
\end{equation}
\end{thm}
\begin{proof} Assume for contradiction that $e>f+g$. Note that
\begin{equation}
\label{aaa}
\frac{g}{f}+\frac{f}{g}-\frac{e^2}{fg}<-2\Longleftrightarrow e>f+g
\end{equation}
by clearing the denominators and extracting full squares.
Similarly
\begin{equation}
\label{bbb}
\frac{e}{f}+\frac{f}{e}-\frac{g^2}{ef}>2\Longleftrightarrow g< |e-f|,
\end{equation}
and both inequalities hold if $e>f+g$.

Denote by $f_1$ and $g_1$ the edges at the other end of $e$.
From above it then follows that
$$0\, <\, eX_e\, =\, \frac{g_1}{f_1}+\frac{f_1}{g_1}-
\frac{e^2}{f_1g_1}+\frac{g}{f}+\frac{f}{g}-\frac{e^2}{fg}\ <\
\frac{g_1}{f_1}+\frac{f_1}{g_1}-\frac{e^2}{f_1g_1}-2$$
using inequality (\ref{aaa}). But this is just the inequality
(\ref{bbb}) for edges $e$, $f_1$ and $g_1$, and thus $e<|f_1-g_1|$.
Suppose $f_1>g_1$; then it follows that $f_1>g_1+e$, and we can apply a
similar argument to the edges $f_2$ and $g_2$ at the other end of $f_1$
to obtain $f_1<|f_2-g_2|$. Continuing this process inductively, and
assuming at each step that $f_i>g_i$, we end up with an infinite
strictly increasing sequence of edges $e<f_1<\ldots <f_n<\ldots$,
which is rather hard to achieve on a finite graph. Thus assuming that
one triangle inequality among $\l$-lengths fails, we have arrived at a
contradiction.
\end{proof}

To demonstrate the power of the triangle inequality we prove a simple 
corollary:

\begin{prop}
\label{rholess}
Let $\r$ be twice the sum of all simplicial coordinates of all edges in
the graph $\Gamma$, as before. Denote the minimal edge in the graph by
$\mu$. If the triangle inequalities are satisfied, then
$\r<\frac{8V}{\mu}$.
\end{prop}
\begin{proof}
For any vertex $v$ denote by $e_v$, $f_v$ and $g_v$ the edges containing
it, with $e_v$ being the maximal among the three. Then using the
triangle inequalities we have
$$\frac{\r}{2}=\sum\limits_v\frac{e_v}{f_vg_v}+
\frac{f_v}{e_vg_v}+\frac{g_v}{e_vf_v}\le\sum\limits_v\frac{f_v+g_v}
{f_vg_v}+\frac{1}{f_v}+\frac{1}{g_v}\le\sum\limits_v \frac{4}{\mu}$$
\end{proof}

\section{Stoke's theorem}
In section 3 we obtained an expression for the \WP form in
terms of $\l$-lengths. However, this expression has multiple summands,
each omitting $n$ variables. In this section we use the Stoke's theorem
to combine the integrals of the summands into one integral of the highest
degree form over a domain in $\R^N$.

\begin{prop}
Let $\w$ be an $(N-n)$-form in $\R^N$. Then
$$\int\limits_{D(\Gamma)} \w=\pm\int\limits_{0\le\r_i\le 1,\ X_e>0}\d\r_1
\wedge\ldots\wedge\d\r_n\wedge\w$$
\end{prop}
\begin{proof}
We apply the Stoke's theorem multiply. Indeed,
$$\int\limits_D \w=\int\limits_{D}\r_1\w=\pm\!\!\!\int
\limits_{X_e>0,\ \r_2=\ldots=\r_n=1,\ 0\le\r_1\le 1}\d\r_1\wedge\w=$$
$$=\!\!\!\int\limits_{X_e>0,\ \r_2=\ldots=\r_n=1,\ 0\le\r_1\le 1}\r_2\d\r_1\wedge\w=
\cdots=\pm\!\!\!\int\limits_{0\le\r_i\le 1,\ X_e>0}\d\r_1\wedge\ldots\wedge\d\r_n
\wedge\w.$$
\end{proof}

Now we need to deal with the $\d\r_i$ factors. We prove the following
\begin{thm}
\label{Stokes}
For $I$ being some set of indices with $|I|=N-n$,
$$\left|\,\int\limits_{0\le\r_i\le 1,\ X_e>0}\d\r_1\wedge\ldots\wedge\d\r_n
\wedge\prod\limits_{i\in I}\frac{\d e_i}{e_i}\right|\le\left|
\,\int\limits_{0\le\r_i\le 1,\ X_e>0}n!\r^n\prod\limits_{i=1}^N
\frac{\d e_i}{e_i}\right|$$
\end{thm}
\begin{proof} Indeed, recall (\cite{wpvol}, section 3.3.4) the definition of the 
$\a$-lengths $\a(e,v):=\frac{e}{fg}$ in the usual notations. Lemma 3.4.2
in \cite{wpvol} states that $\rho_i$ is twice the sum of $\a$-lengths of
the sectors it traverses. What matters for us is that it is a sum of
some $\a$-lengths. For any $\a$-length we have
$$\frac{\p\a(e,v)}{\p e}=\frac{1}{e}\a(e,v)\quad{\rm and}\quad
\frac{\p\a(e,v)}{\p f}=-\frac{1}{f}\a(e,v).$$
Thus for $\r_i=2\sum\limits_{j=1}^n\alpha(f_i,v_i)$ we have
$$\left|\frac{\p\r_i}{\p e}\right|=
=\frac{2}{e}\left|\sum\pm\a(f_j,v_j)\right|\le\frac{2}{e}\sum\a(f_j,v_j)
=\frac{\r_i}{e}<\frac{\r}{e}.$$
Applying this trick to each of $\d\r_i$ yields the theorem.
\end{proof}
Combining this theorem with the bound on the coefficients of the \WP from 
theorem \ref{aI},  noting that there are ${N\choose n}<N^n/n!$
summands in the \WP form, and enlarging $D$ to the domain 
$0\le\r\le n, X_e>0\ \forall e$, we get the following
\begin{cor}
\label{form}
The integral of the \WP volume over the domain of integration is bounded by
$$\int\limits_D\w_{WP}^{\wedge(3g-3+n)}<2^NN^n\left|\,\int\limits_{
0\le\r<n,\ X_e>0}\r^n\prod\limits_{i=1}^N\frac{\d e_i}{e_i}\right|$$
\end{cor}

\section{Triangle inequality combinatorics}
Now we proceed to show how the triangle inequalities lead to a
converging integral as an upper bound for the \WP volume. 
We develop an algorithm for inductive estimation of the integral of the \WP
volume form over the domain where the triangle inequalities hold. 

\begin{df} Two edges $e$ and $f$ of the graph are called {\it linked} if
they are adjacent at a vertex $v$, and the third edge $g$ at $v$ is the
minimal of the three (not necessarily strictly). Notice that from the
triangle inequalities $e<f+g$ and $f<e+g$ it then follows that
$\frac{1}{2}\le \frac{e}{f}\le 2$.
\end{df}

\begin{df} A {\it chain} is a sequence of edges in which every
two consecutive ones are linked. We will only be interested in
maximal chains --- the ones which are not a part of any longer chain. Such a
chain must either form a loop in the graph, or end by two edges which are
minimal at their outer-end vertices.
\end{df}

\begin{df} Define a {\it wheel} to be an ordered sequence of (maximal) chains
$c_1,\ldots, c_m$ such that for any $i$ there is at least one chain
among $c_1,\ldots, c_{i-1}$ ending at a vertex inside the chain $c_i$. In
further considerations, we will only be interested in maximal wheels,
the ones which cannot be enlarged any further. Basically this means that
the ends of all chains in the wheel belong to other chains
already included in the wheel, so that there are no edges ``sticking
out'' of the wheel.
\end{df}

For further computations we split the domain of
integration into at most $3^V$ parts by deciding which two edges
at each vertex are linked. We then use chains and wheels to introduce some order on the set
of edges, and to integrate inductively. In this we will be aided by the
following technical lemmas. The notations are as in the definition of linking: $e$,
$f$ and $g$ are the three edges at some vertex, among which $g$ is
minimal. We are working in the domain $\Delta$ where the triangle inequalities
hold, noting that $\Delta\supset D(\Gamma)$ by theorem \ref{le2thm}. 

\begin{lm}
\label{trick1}
For fixed $e$ we have $\int_\Delta \frac{\d f}{f}\le\ln 4$,
where $\int_\Delta$ denotes integration over the possible values of $f$ within
domain $\Delta$ for fixed $e$.
\end{lm}
\begin{proof}
Since $e$ and $f$ link, we have $e/2<f<2e$. Thus
$$\int_\Delta \frac{\d f}{f}\le\int_{e/2}^{2e}\frac{\d f}{f}=\ln 4.$$
\end{proof}

\begin{cor}
\label{chainok}
Let $e_1\ldots e_m$ be a chain. If we fix (the $\l$-length of) an edge
$e_i$, then
$$\int_\Delta\prod\limits_{j\ne i}\frac{de_j}{e_j}\le(\ln 4)^{m-1}$$
\end{cor}
\begin{proof}
Start from the ends of the chain, and apply the lemma to eliminate edges
one by one, coming from the ends towards $e_i$.
\end{proof}

\begin{lm}
\label{trick2}
If we fix the edge $g$, then the integral over the linked edges $e$ and
$f$ can be estimated as $\int_\Delta\frac{\d e\d f}{ef}<2$
\end{lm}
\begin{proof}
Split the integral into two parts depending on whether $e>f$ or
$f>e$. The computation for them is identical; if $f<e$, we have
$$\int\limits_{\Delta\cap\lbrace f<e\rbrace}\frac{\d e\d f}{ef}\le
\int\limits_g^\infty\frac{\d f}{f}\int\limits_f^{f+g}\frac{\d
e}{e}=\int\limits_g^\infty\frac{\d f}{f}\ln\left(1+\frac{g}{f}
\right)<\int\limits_g^\infty\frac{\d f}{f}\frac{g}{f}=1.$$
\end{proof}

\begin{prop}
\label{wheelok}
Starting from two linked edges $e$ and $f$, construct a wheel consisting 
of edges $e_1\ldots e_m$, where the list includes $e$ and $f$ themselves. 
Then for fixed $g$
$$\int_\Delta\prod\limits_{i=1}^m\frac{\d e_i}{e_i}\le (\max(\ln
4,\sqrt 2))^m$$
\end{prop}
\begin{proof}
The wheel is a collection of chains $c_1\ldots c_k$. Keep the two edges
of $c_k$ in between which some chain $c_i$ ends (existent by
definition of a wheel), and integrate the other ones out using
corollary \ref{chainok} --- by this we pick up a factor of $\ln 4$
for each edge. Then use lemma \ref{trick2} above to integrate out
the last two remaining edges of the chain $c_k$ --- here we pick
up a factor of $2$ for two edges, i.e.~$\sqrt 2$ per edge.
Performing induction in $k$ finishes the proof.
\end{proof}

If a wheel were the whole graph, we would be able to estimate the
integral using the above proposition. However, if the wheel is not the
whole graph, we need to be able to link it to the rest of the graph.
Thus we will need the following

\begin{lm}
\label{trick3}
In the usual notations for fixed $e$ we have
$$\int_\Delta \frac{\d f\d g}{fg}<\frac{8}{3}$$
\end{lm}
\begin{proof}
Indeed, recall that $g>4$ by proposition \ref{ge4}. Thus
$$\int\limits_4^{e}\frac{\d g}{g}\int\limits_{\max(g,e-g)}^g \frac{\d
f}{f}+\int\limits_4^g\frac{\d g}{g}\int\limits_e^{e+g}\frac{\d
f}{f}\le\int\limits_4^{e/2}\frac{\d g}{g}\ln \frac{e}{e-
g}+\int\limits_4^e\frac{\d g}{g}\ln \left(1+\frac{g}{e}\right)$$
Since $\ln(1+x)<x$ for $x>0$,  for the second summand we have
$$\int\limits_4^e\frac{\d g}{g}\ln(1+\frac{g}{e})<\int\limits_4^e
\frac{\d g}{g}\frac{g}{e}=\frac{e-4}{e}<1.$$
For the first summand we compute
$$\int\limits_4^e \frac{\d g}{g}\ln\frac{e}{e-g}=-\int\limits_4^e
\frac{\d g}{g}\ln\left(1-\frac{g}{e}\right)=\int\limits_4^e\frac{\d g}{g}
\sum\limits_{n=1}^\infty\frac {g^n}{n e^n}=$$ $$=\sum\limits_{n=1}^\infty
\frac{1}{n e^n}\int\limits_4^e g^{n-1}\d g=\sum\limits_{n=1}^\infty
\frac{1}{n e^n}\frac{e^n-4^n}{n}<\sum\limits_{n=1}^\infty\frac{1}{n^2}<
\frac{5}{3}.$$
Combining the above estimates, we get the lemma.
\end{proof}

\begin{thm}
\label{intout}
Let $\mu$ be the minimal edge of $\Gamma$, and let $e_1\ldots e_{N-1}$ be
all the other edges of the graph. Then for a fixed value of $\mu$ and a 
fixed choice of the two linking edges at each vertex we have
$$\int_\Delta\prod\limits_{i=1}^{N-1}\frac{\d
e_i}{e_i}<\left(\frac{8}{3}\right)^{(N-1)/2}$$
\end{thm}
\begin{proof}
Construct a wheel $w_1$ starting from edge $\mu$. If this wheel is
not the whole graph, consider a chain $c_{1,1}$ ending on $w_1$ by at
least one end. If it ends on $w_1$ by the other end also, consider
another such chain $c_{1,2}$ and so on, until we either exhaust
the graph, or get a chain $c_{1,m_1}$ which has an end not on $w_1$.
Then construct a wheel $w_2$ at the other end of the chain
$c_{1,m_1}$. If the union of these two wheels and the chains
constructed is still not the whole graph, we repeat the process.

As a result, we decompose the graph into a disjoint union of wheels
$w_1\ldots w_k$ and chains $c_{i,j}$ for $i\le k$ and $j\le m_i$
such that the chains $c_{i,j}$ have both ends on $w_i$ for $i=k$ or 
$j<m_i$, and that $c_{i,m_i}$ connects $w_i$ and $w_{i+1}$ for $i<k$.
Then use corollary \ref{chainok}
to eliminate all edges of chains $c_{k,i}$ except the terminal
ones, which end at $w_k$. Using lemma \ref{trick3}, we can
then include these terminal edges of $c_{k,i}$'s while
implementing the proof of proposition \ref{wheelok} --- integrating out 
the edges of $w_k$ one by one. Doing this, we
get a factor of $8/3$ for eliminating two edges, instead of the 
smaller factor of $(\ln 4)^2$, which we were getting originally. 

At the last step
of integrating over the edges of $w_k$ we use the edge at the end
of $c_{k-1,m_{k-1}}$ for lemma \ref{trick2}, and thus reduce the
problem to $k-1$ wheels. Induction in $k$ then yields the desired
result.
\end{proof}

Now we combine all the above estimates to finally obtain

\begin{thm}
\label{result}
In the above notations,
$$\int_{D(\Gamma)}\w_{WP}^{\wedge(3g-3+n)}<2^N3^VN^n\left(\frac{8}{3}\right)^{(N-1)/2}(2V)^n$$
\end{thm}
\begin{proof}
Combining the results of corollary \ref{form} and proposition \ref{rholess},
we see that the integral in question is bounded above by
$$2^NN^n\left|\,\int \r^n\frac{\d\mu}{\mu}
\prod\limits_{i=1}^{N-1}\frac{\d e_i}{e_i}\right|<
2^NN^n\left|\,\int \left(\frac{8V}{\mu}\right)^n\frac{\d\mu}{\mu}
\prod\limits_{i=1}^{N-1}\frac{\d e_i}{e_i}\right|.$$
Using theorem \ref{intout}, we can integrate out 
all variables except $\mu$, by acquiring an extra factor of $(8/3)^{(N-1)/2}$. 
Remembering the factor of $3^V$ for choosing the minimal edge at each vertex,
our final upper bound becomes
$$2^N3^VN^n\left(\frac{8}{3}\right)^{(N-1)/2}(8V)^n\int\limits_4^\infty
\frac{\d\mu}{\mu^{n+1}}=2^N3^VN^n\left(\frac{8}{3}\right)^{(N-1)/2}(2V)^n$$
\end{proof}

Using our explicit knowledge of the \WP volume form, in the case of one 
puncture we get
$$\int\limits_D \w_{WP}^{\wedge(3g-3+n)}<2^{4g-2} 3^V(8/3)^{(N-1)/2} N.$$

\section{Conclusion}
For the case of one puncture, combining our estimates with Penner's 
asymptotic computation of the number of
cells in formula \ref{domain} (being $\frac {(2g)!}{N}\left(e/6\right)^{-2g}$), and estimating
$1/{\rm Aut}\,\Gamma$  from above by one, we finally get
$${\rm vol}_{WP}(\M_{g,1})< (2g)!\, 2^{4g-2}3^{4g-2} (\ln 4)^{6g-3}
\left(\frac{6}{e}\right)^{2g}=:c_1^g(2g)!,$$
which has the same leading order infinity as Penner's lower bound
$${\rm vol}_{WP}(\M_{g,1})>\left(\frac{8e^2}{9}\right)^{2g}\frac{(2g)!}
{2(6g-3)^2}=:c_2^g(2g)!,$$ 
where $c_i$ are some explicit constants. Thus we have proven that
$$c_2^g(2g)!<{\rm vol}_{WP}(\M_{g,1})<c_1^g(2g)!\quad{\rm and}\quad
\lim\limits_{g\to\infty}\frac{\ln{\rm vol}_{WP}(\M_{g,1})}{g\ln g}=2.$$

Intuitively, for more than one puncture the number of graphs should grow 
with genus in the same way as for one puncture, since all the punctures are 
far away from the additional handles being added, and their number should not
matter. 

Rigorously, let $T(g,n)$ denote the set of isomorphism classes of 
ideal triangulations of a surface of genus $g$ with $n$ punctures. 
Then we prove
\begin{prop}
There is a following upper bound on the number of triangulations:
$$|T(g,n)|<\frac{N^2}{2}|T(g,n-1)<\ldots<\frac{N^{2n-2
}}{2^{n-1}}|T(g,1)|<\frac{(2g)!N^{2n-3}}{2^{n-1}}
\left(\frac{6}{e}\right)^{2g}$$
\end{prop}
\begin{proof}
For $n>1$  we construct a relation $\phi\in T(g,n)\times T(g,n-1)$
in the following way.  Consider two distinct punctures $p_1$ and $p_2$
connected by an edge $e$ of a triangulation $x\in T(g,n)$ --- if 
such did not exist, i.e. if all edges emanating from a puncture went back 
to the puncture itself, it would not be a triangulation of the surface. 
Shrinking $e$ to a point, and collapsing triangles on both sides of $e$ into 
arcs of a triangulations, thus identifying $p_1$ and $p_2$, produces a new
triangulation $y\in T(g,n-1)$. We define $\phi$ to be the set of all pairs
$(x,y)$ obtained in such a way.

Now consider some $y\in T(g,n-1)$. Any graph in $x\in T(g,n)$ such that
$(x,y)\in\phi$ can be reconstructed from $y$ by ``blowing up'' a pair
of edges emanating from a vertex to triangles. Since there are 
$(6g-6+3(n-1))(6g-6+3(n-1)-1)/2<N^2/2$ ways to choose a pair of edges of $y$,
there are at most $N^2/2$ triangulations $x\in T(g,n)$ such that $(x,y)\in\phi$. The argument works for all $y$, and thus
$|T(g,n)|<N^2|T(g,n-1)|/2$. Applying this argument until we decrease the 
number of punctures to one, and then utilizing Penner's asymptotic computation
for that case finishes the proof. \end{proof}

Combining all our estimates, for $g\gg n$ we get
$${\rm vol}_{WP}(\M_{g,n})<\frac{(2g)!N^{2n-3}}{2^{n-1}}\left(\frac{6}{e} 
\right)^{2g}2^N3^VN^n\left(\frac{8}{3}\right)^{N/2}(2V)^n<C^g(2g)!,$$
where $c$ is any constant greater than $2^{17}3^3/e^2$ (notice that is is 
independent of $n$, since the $g^{n}$ has a lower growth order), and thus
$$\lim\limits_{g\to\infty,\ n \ {\rm fixed}}\frac{\ln {\rm vol}_{WP}
(\M_{g,n})}{g\ln g}\le 2.$$

\section*{Acknowledgements}
The author would like to thank Professor Yum-Tong Siu for suggesting the 
problem and for many invaluable discussions. Without Professor Siu's continuous
support this work would have been impossible. We would also like to thank the 
University of Hong Kong for hospitality during November 1999, when this work 
was finalized.


\begin{thebibliography}{9}
\bibitem{fabnum}{C.~Faber {\it Algorithms for computing intersection numbers
on moduli spaces of curves, with an application to the class of the locus of
Jacobians}. New trends in algebraic geometry (Warwick, 1996), 93--109, London
Math. Soc. Lecture Note Ser., 264, Cambridge Univ. Press, Cambridge,
1999.}
\bibitem{kontsevich}{M.~Kontsevich {\it Intersection theory on the
moduli space of curves} (Russian) Funktsional. Anal. i
Prilozhen. {\bf 25} (1991), no. 2, 50--57; translation in Functional
Anal. Appl. {\bf 25} (1991), no. 2, 123--129}
\bibitem{mz}{Yu.~I.~Manin, P.~Zograf {\it Invertible Cohomological Field
Theories and Weil-Petersson volumes}, preprint math.AG/9902051 v2}
\bibitem{wpvol}{R.~C.~Penner {\it \WP volumes}, J.~Differential Geometry
{\bf 35} (1992) 559-608}
\bibitem{intscheme}{R.~Penner {\it The decorated Teichm\"uller space
of punctured surfaces}, Comm. Math. Phys. {\bf 113} (1987) 299--319}
\bibitem{witten}{E.~Witten {\it Two-dimensional gravity and intersection
theory on moduli space}, Surveys in Differential Geometry {\bf 1}
(1991), 243--310}
\bibitem{wolpert}{S.~Wolpert {\it On the homology of the moduli space of
stable curves}, Annals of Mathematics {\bf 118} (1983), 491--523}
\bibitem{zograf}{P.~Zograf {\it \WP volumes of moduli spaces of curves
and the genus expansion in two dimensional gravity}, preprint
math.AG/9811026}
\end{thebibliography}
\end{document}